\documentclass{article}

\usepackage{amsmath}
\usepackage{amsfonts}
\usepackage{amssymb}
\usepackage{amsthm}

\newcommand{\N}{\mathbb{N}}

\newcommand{\bs}{\boldsymbol}

\newcommand{\E}{\mathbb{E}}

\renewcommand{\Pr}{\mathbb{P}}

\theoremstyle{definition}

\newtheorem{proposition}{Proposition}

\theoremstyle{remark}

\begin{document}

\title{Upward and Downward Runs on Partially Ordered Sets}

\author{Kyle Siegrist\\
Department of Mathematical Sciences \\
University of Alabama in Huntsville}

\maketitle

\begin{abstract}
We consider Markov chains on partially ordered sets that generalize the success-runs and remaining life chains in reliability theory.  We find conditions for recurrence and transience and give simple expressions for the invariant distributions.  We study a number of special cases, including rooted trees, uniform posets, and posets associated with positive semigroups.

\vspace{5 pt}

\vspace{5 pt}

\end{abstract}

\section{Partially Ordered Sets}

\subsection{Preliminaries}

Suppose that $(S, \preceq)$ is a discrete partially ordered set.  Recall that $C \subseteq S$ is a {\em chain} if $C$ is totally ordered under $\preceq$.  We make the following assumptions:
\begin{enumerate}
\item \label{a.m} There is a minimum element $e$.
\item \label{a.f} For every $x \in S$, every chain in $S$ from $e$ to $x$ is finite.
\end{enumerate}

Recall that $y$ {\em covers} $x$ if $y$ is a minimal element of $\{t \in S: t \succ x\}$. The {\em covering graph} (or Hasse graph) of $(S, \preceq)$ is the directed graph with vertex set $S$ and edge set $E = \{(x, y) \in S^2: y \text{ covers } x\}$. From the assumptions, it follows that for each $x \in S$, there is a (directed) path from $e$ to $x$ in the graph, and every such path is finite. For $x \in S$, let
\[ A_x = \{y \in S: y \text{ covers } x\}, \quad B_x = \{w \in S: x \text{ covers } w\} \]
That is, $A_x$ is the set of elements immediately {\em after} $x$ in the partial order, while $B_x$ is the set of elements immediately {\em before} $x$ in the partial order. Note that $A_x$ could be empty or infinite. On the other hand, $B_e = \emptyset$, but for $x \ne e$, $B_x \ne \emptyset$ since there is at least one path from $e$ to $x$.

An {\em upward run chain} on $(S, \preceq)$ is a Markov chain that, at each transition, moves to a state immediately above the current state or back to $e$, the minimum state.  A {\em downward run chain} is a Markov chain that, at each transition, moves to a state immediately below the current state, unless the current state is $e$ in which case the chain can move anywhere in $S$.  For particular posets, upward and downward runs can have applications in reliability theory, communications theory, queuing theory and other areas.  Generally, posets are the  natural mathematical home for these stochastic processes. In this article, we are interested in general issues of recurrence, invariant distributions, time reversal, and results for special types of posets. See \cite{Evstigneev} for another class of Markov chains on posets.

\subsection{Uniform posets} \label{ss.uniform0}

An interesting case is when the partially ordered set $(S, \preceq)$ is {\em uniform} in the sense that for each $x \in S$, all paths from $e$ to $x$ have the same length. It then follows that the all paths from $x$ to $y$ have the same length for any $x, \, y \in S$ with $x \preceq y$; we denote this length by $d(x, y)$. Let $S_n = \{x \in S: d(e, x) = n\}$ for $n \in \N$.  Of course, $S_0 = \{e\}$ and $\{S_n: n \in \N\}$ partitions $S$. 

\subsection{Rooted trees and path space} \label{ss.paths}

Another important special case is when the covering graph of $(S, \preceq)$ is a rooted tree with root $e$.  In this case, $B_x$ has a single element, which we will denote by $x^-$, for each $x \ne e$. There is a unique path from $e$ to $x$ for each $x \in S$. Thus the poset $(S, \preceq)$ is uniform, so the definitions in Section \ref{ss.uniform0} apply. 

In fact, rooted trees form an essential special case, because we will show that upward and downward runs on an arbitrary poset can be constructed from upward and downward runs on a certain rooted tree of paths. Specifically, suppose that $(S, \preceq)$ is a poset, and let $\hat S$ denote the set of finite, directed paths in $S$, starting at $e$.  We define the partial order $\preceq$ on $\hat S$ by $a \preceq b$ if and only if $a$ is a prefix of $b$. The covering graph of $(\hat S, \preceq)$ is a tree rooted at $e$ (the degenerate path consisting only of $e$).  For $a \in \hat S$, let $m(a)$ denote the endpoint of $a$.  If $y$ covers $x$ in $S$ then for every $a \in \hat S$ with $m(a) = x$, $ay$ covers $a$ (where $ay$ denotes that path obtained by appending $y$ to the end of $a$). Let $\hat S(x) = \{a \in \hat S: m(a) = x\}$, the set of paths ending in $x$. Note that for $a \ne e$, $a^-$ is the path obtained by removing the endpoint of $a$.  If $a \in \hat S$, the unique path from $e$ to $a$ in $\hat S$ simply consists of the successive prefixes of $a$.  Thus in the notation of Section \ref{ss.uniform0}, $d(e, a)$ is the length of the path $a$ and hence $\hat S_n$ is the set of paths of length $n$.

\subsection{Positive semigroups} \label{ss.semigroups0}

Another important special case is when the partially ordered set $(S, \preceq)$ is associated with a {\em positive semigroup} $(S, \cdot)$.  That is, $\cdot$ is an associative binary operation on $S$ with an identity element $e$, no non-trivial inverses, and satisfying the left-cancellation law. In this case, $x \preceq y$ if and only if there exists $t \in S$ with $xt = y$. Positive semigroups are essentially characterized by the fact that $xS = \{y \in S: y \succeq x\}$ is order-isomorphic to $S$ for each $x \in S$; the mapping $t \mapsto xt$ is an isomorphism. Assumption \ref{a.m} is always satisfied since $e$ is the minimum element.  We will assume that $[e, x] = \{t \in S: t \preceq x\}$ is finite for each $x \in S$, so Assumption \ref{a.f} is satisfied as well. Probability distributions on positive semigroups are studied in \cite{Rowell, Siegrist, Siegrist2, Siegrist3, Siegrist4}

An element $i \in S$ is {\em irreducible} if $i$ cannot be factored, except for the trivial factoring $i = ie = ei$. If $I$ is the set of irreducible elements of $(S, \cdot)$ then $A_x = \{xi: i \in I\}$, so in particular, $\#(A_x)$ is the same for each $x$. The poset $(S, \preceq)$ will be uniform if and only if for every $x \in S$, all factorings of $x$ over $I$ have the same number of factors.  In this case, $d(e, x)$ is the number of factors.

The path space associated with $(S, \cdot)$ is isomorphic to the {\em free semigroup} on the set of irreducible elements $I$. This is the set of finite ``words'' over the alphabet $I$, with concatenation as the semigroup operation. That is, a path from $e$ to $x$ in path space is uniquely associated with a factoring of $x$ over $I$: $x = i_1 i_2 \cdots i_n$; the string on the right is a word in the free semigroup. 

\subsection{Probability Distributions}

Suppose that $(S, \preceq)$ is a poset and that $X$ is a random variable with support $S$. As usual, the probability density function (PDF) of $X$ is the function $f$ given by $f(x) = \Pr(X = x)$ for $x \in S$. The {\em upper probability function} (UPF) of $X$ is the function $F$ given by 
\[ F(x) = \Pr(X \succeq x) = \sum_{y \succeq x} f(y), \quad x \in S \]
Finally, the {\em rate function} of $X$ is the function $r$ given by 
\[ r(x) = \frac{f(x)}{F(x)} = \Pr(X = x | X \succeq x), \quad x \in S \]
In particular, $X$ has {\em constant rate} if $r$ is constant on $S$. For general posets, the distribution of $X$ is not uniquely determined by the UPF $F$ (and certainly not by the rate function $r$). These issues and the existence of constant rate distributions are explored in \cite{Siegrist5}. A special case of a general expected value result in \cite{Siegrist5} is
\begin{equation} \label{eq.expect}
\sum_{x \in S} F(x) = \E[\#(D[x])]
\end{equation}
where $D[x] = \{t \in x: t \preceq x\}$. If $(S, \preceq)$ is a rooted tree, then the UPF $F$ of $X$ does determine the distribution of $X$, since clearly the PDF $f$ of $X$ is given by
\[ f(x) = F(x) - \sum_{y \in A(x)} F(y) \]
Moreover, when $(S, \preceq)$ is a rooted tree, $\#(D[x]) = 1 + d(e, X)$

\section{The Upward Run Chain} \label{s.upward}

\subsection{Basic definitions and results}

A Markov chain $\bs{X} = (X_0, X_1, \ldots)$ on a poset $(S, \preceq)$ is an {\em upward run} chain if the transition function $P$ satisfies $P(x, y) > 0$ if and only if $y \in A_x$ or  $y = e$.  Of course, we must have
\[ P(x, e) = 1 - P(x, A_x), \quad x \in S\]
Thus, in state $x \in S$, the chain next moves to a state $y \in A_x$, or back down to $e$. The chain is irreducible since $e$ leads to every state, and every state leads back to $e$. The chain is aperiodic since $e$ leads back to $e$ in one step. If the state space is $(\N, \le)$, then an upward run chain is simply a {\em success-runs} chain; these are commonly studied in reliability theory.

For $x \in S$, let $T_x$ denote the first positive hitting time to $x$:
\[ T_x = \min \{n \in \N_+: X_n = x\} \]
Define $F: S \to (0, 1]$ by $F(x) = \Pr_e(T_x \le T_e)$. Of course, $\Pr_e(T_x = T_e) = 0$ if $x \ne e$. Thus, $F(x) = \Pr_e(T_x < T_e)$ if $x \ne e$ while $F(e) = 1$. 

\begin{proposition} 
Suppose that $\bs{X}$ is recurrent. Then $F$ is left-invariant for $P$.
\end{proposition}

\begin{proof}
Suppose first that $y \ne e$. Starting at $e$, the chain moves upward to $y$ without an intermediate return to $e$ if and only if the chain moves upward to some $x \in B_y$, without an intermediate return to $e$, and then moves in one step from $x$ to $y$. Therefore
\[ F(y) = \sum_{x \in B(y)} F(x) P(x, y) = (F P)(y) \]
On the other hand, to return to $e$, starting at $e$ the chain must go directly back to $e$ or move upward to some intermediate state $x \in S$ and then go back to $e$ in one step. Thus, since the chain is recurrent,
\[ (F P)(e) = \sum_{x \in S} F(x) P(x, e) = \Pr_e(T_e < \infty) = 1 = F(e) \]
\end{proof}

We will refer to $F$ as the {\em standard invariant function} for $\bs{X}$. In the recurrent case, the chain is positive recurrent if and only if
\begin{equation} \label{eq.precurrent}
\mu(e) := \sum_{x \in S} F(x)
\end{equation}
is finite. In this case, the invariant PDF $f$ is given by $f(x) = F(x) / \mu(e)$ for $x \in S$. Equivalently, $\mu(x) := \E_x(T_x) = \mu(e) / F(x)$ for $x \in S$.

\subsection{Upward runs on rooted trees} \label{ss.upTrees}

Consider now the special case where the covering graph of $(S, \preceq)$ is a rooted tree.  If $x \in S$ and $e \, x_1 \, x_2 \cdots x_{n-1} \, x$ is the unique path in $S$ from $e$ to $x$, then clearly
\begin{equation} \label{eq.up1}
F(x) = P(e, x_1) P(x_1, x_2) \cdots P(x_{n-1},x)
\end{equation}
the product of the transition probabilities along the path from $e$ to $x$.
It follows that
\begin{equation} \label{eq.up2}
\Pr_e (T_e > n) = \sum_{x \in S_n} F(x), \quad n \in \N
\end{equation}
Thus the chain is recurrent if and only if the sum in (\ref{eq.up2}) has limit 0 as $n \to \infty$.

Suppose that the chain is recurrent.  Then clearly $F$ is the UPF of $X_{T(e) - 1}$, the last state vistied before returning to $e$ (starting at $e$). That is
\[ F(x) = \Pr_e(T_x \le T_e) = \Pr_e(X_{T(e) - 1} \succeq x), \quad x \in S \]
Note from (\ref{eq.expect}) that $\mu(e) = \E_e(T_e)$ can be written as 
\[ \mu(e) = 1 + \E_e \left[ d \left( e, X_{T(e) - 1} \right) \right] \]
In the positive recurrent case, the invariant PDF $f$ is the function obtained by normalizing $F$ with $\mu(e)$. Note that in general, $f$ is {\em not} the PDF of $X_{T(e) - 1)}$ given $X_0 = e$. In fact, the two PDFs are the same if and only if the invariant distribution has constant rate.

Conversely, given an UPF $F$ on $S$, we can construct a recurrent upward run chain with $F$ as the standard invariant function. Specifically, suppose that $X$ is random variable with support $S$ and with UPF $F$ and PDF $f$. Define $P$ by
\begin{align*}
	P(x, y) &= \Pr(X \succeq y | X \succeq x) = \frac{F(y)}{F(x)}, \quad x \in S, \, y \in A_x \\
	P(x, e) &= \Pr(X = x | X \succeq x) = \frac{f(x)}{F(x)}, \quad x \in S
\end{align*}
Then $P$ is the transition probability function for an upward run chain $\bs{X} = (X_0, X_1, \ldots)$. Moreover, $F$ and $P$ satisfy (\ref{eq.up1}), so the notation is consistent and $F$ is the standard invariant function.  Note that $x \mapsto P(x, e)$ is the rate function of $X$. In particular, if $P(x, e) = \alpha$ for all $x \in S$, then $X$ has constant rate $\alpha$, the chain is positive recurrent, and $f$ is the invariant PDF.

\subsection{Upward runs on path space} 

Now return to the general case where $(S, \preceq)$ is an arbitrary poset and $(\hat S, \preceq)$ is the corresponding path space discussed in Section \ref{ss.paths} (a rooted tree).  Suppose that $P$ is the transition probability function of an upward run chain on $S$. We define  $\hat P$ on $\hat S$ by
\[ \hat P(a, ay) = P(x, y), \; \hat P(a, e) = P(x, e), \quad a \in \hat S(x), \, x \in S \]
Clearly $\hat P$ is a valid transition probability function and corresponds to an upward run chain on $\hat S$.  Moreover, starting at $e$, we can define the two chains on a common probability space.  To do this, we need some notation. If $\bs{\hat X}$ is a process on $\hat S$, define the process $\bs{X}$ on $S$ by $X_n = m(\hat X_n)$.

\begin{proposition} \label{p.upCouple}
If $\bs{\hat X}$ is an upward run chain on $\hat S$ with transition probability function $\hat P$ and starting at $e$, then $\bs{X}$ is an upward run chain on $S$ with transition probability function $P$ and starting at $e$.
\end{proposition}

\begin{proof}
The proof follows easily from the definition of the path space $\hat S$ and the fact that
\[ \hat P(a, b) = P(m(a), m(b)), \quad a, \, b \in \hat S \]
\end{proof}

Of course the results of Section \ref{ss.upTrees} apply to $\bs{\hat X}$. Also, with the chains coupled as in Proposition \ref{p.upCouple}, note that $T_e = \hat T_e$. If $a = e \, x_1 \, x_2 \cdots x_n \in \hat S$ then from (\ref{eq.up2}) and the definition of $\hat P$,
\begin{align*}
\hat F(a) &= \hat P(e, e x_1) \hat P(e x_1, e x_1 x_2) \cdots \hat P(e x_1 \cdots x_{n-1}, e  x_1 \cdots x_n)\\
&= P(e,x_1) P(x_1, x_2) \cdots P(x_{n-1}, x_n)
\end{align*}
Now, for $\bs{X}$ to go from $e$ to $x$ without an intermediate return to $e$, $\bs{X}$ must move along some path from $e$ to $x$ so
\[ F(x) = \sum_{a \in \hat S(x)} \hat F(a), \quad x \in S \]
It then follows that
\begin{equation} \label{eq.path2}
\Pr_e (T_e > n) = \sum_{a \in \hat S_n} \hat F(a)
\end{equation}
The chains $\bs{X}$ and $\bs{\hat X}$ have the same classification: recurrent if and only if the sum in (\ref{eq.path2}) has limit $0$ as $n \to \infty$, and of course positive recurrent if and only if the sum in (\ref{eq.precurrent}) is finite. In the positive recurrent case, the invariant PDFs $\hat f$ and $f$ are related by
\[ f(x) = \sum_{a \in \hat S(x)} \hat f(a), \quad x \in S \]

Since $(\hat S, \preceq)$ is a rooted tree, $\hat F$ is an UPF if the chains are recurrent.  On the other hand $F$ may not be an UPF on $S$, and conversely, there may exist UPFs $F$ on $S$ that cannot be realized as the standard invariant function for a recurrent upward run chain on $S$.

\section{The Downward Run Chain} \label{s.downward}

\subsection{Basic definitions and results}

A Markov chain $\bs{Y} = (Y_0, Y_1, \cdots)$ on $(S \preceq)$, is a {\em downward run chain} if the transition probability function $Q$ satisfies $Q(x, y) > 0$ if and only if $x = e$ and $y \in S$, or if $x \in S - \{e\}$ and $y \in B_x$. Thus, the chain moves downward from a state to one of its predecessors, until it gets to $e$; then it can move anywhere in $S$.  The chain is irreducible since every state leads to $e$ and $e$ leads to every state.  The chain is aperiodic since $e$ leads to $e$ in one step. If the poset is $(\N, \le)$, then a downward run chain is simply a {\em remaining life} chain; these are commonly studied in reliability theory.

We denote the first positive hitting time to $x \in S$ by
\[ U_x = \min \{n \in \N_+: Y_n = x\} \]
A downward chain is always recurrent.  Since all paths from $e$ to $x$ are finite, $\Pr_x(U_e < \infty) = 1$ for $x \in S - \{e\}$. Hence
\[ \Pr_e(U_e < \infty) = Q(e, e) + \sum_{x \in S - \{e\}} Q(e, x) \Pr_x (U_e < \infty) = 1 \]
Define $G: S \to (0, 1]$ by $G(x) = \Pr_e(U_x \le U_e)$. Of course, $\Pr_e(U_x = U_e) = 0$ if $x \ne e$. Thus, $G(x) = \Pr_e(U_x < U_e)$ if $x \ne e$ while $G(e) = 1$. Note that $G$ is defined for the downward run chain just as $F$ is defined for the upward run chain.

\begin{proposition} 
The function $G$ is left-invariant for $Q$.
\end{proposition}

\begin{proof}
For $y \in S$, 
\begin{align*}
(G Q)(y) &= \sum_{x \in S} G(x) Q(x, y)\\
&= G(e) Q(e, y) + \sum_{x \in A(y)} G(x) Q(x, y)\\
&= Q(e, y) + \sum_{x \in A(y)} G(x) Q(x, y) = G(y)
\end{align*}
For the last line, note that starting at $e$, the chain hits $y$ before returning to $e$ if and only if the chain jumps immediately to $y$ or hits some $x \in A_y$, before returning to $e$, and then moves from $x$ to $y$ in one step.
\end{proof}

As before, we will refer to $G$ as the {\em standard invariant function} for $\bs{Y}$. The chain is positive recurrent if and only if 
\[\nu(e) := \sum_{x \in S} G(x)\]
is finite.  In this case, the invariant PDF is given by $g(x) = G(x) / \nu(e)$ for $x \in S$. Equivalently, $\nu(x) := \E_x(U_x) = \nu(e) / G(x)$ for $x \in S$. 

\subsection{Downward runs on rooted trees} \label{ss.downTrees}

Consider the special case where the covering graph of $(S, \preceq)$ is a rooted tree.  Since a non-root vertex has a single parent, we have
\[ Q(x, x^-) = 1, \quad x \in S - \{e\} \]
Thus, the invariant function $G$ is simply the UPF of $Y_1$, given $Y_0 = e$:
\[ G(x) = \sum_{y \succeq x} Q(e, y) = \Pr_e(Y_1 \succeq x), \quad x \in S \]
From (\ref{eq.expect}), $\nu(e) = \E_e(U_e)$ can be written as
\[\nu(e) = 1 + \E_e [d(e, Y_1)] \]
Note again that in the positive recurrent case, the invariant PDF is the function $g$ obtained by normalizing $G$, and is not in general the PDF of $Y_1$ given $Y_0 = e$. The two PDFs are the same if and only if the distribution has constant rate.

Conversely, given an UPF $G$ on $S$, it's trivial to construct a downward run chain $\bs{Y}$ with $G$ as the standard invariant function.  Specifically, suppose that $G$ is the UPF corresponding to the PDF $g$. We just need to define $Q(x, e) = g(x)$ for $x \in S$, and of course $Q(x, x^-) = 1$ for $x \in S - \{e\}$.

\subsection{Downward runs on path space} 

Now return to the general case where $(S, \preceq)$ is an arbitrary poset and $(\hat S, \preceq)$ is the corresponding path space discussed in Section \ref{ss.paths} (a rooted tree).  Suppose that $Q$ is the transition probability function of a downward run chain on $S$. We define  $\hat Q$ on $\hat S$ as follows: 
\begin{align*}
	\hat Q(e, e \, x_1 \cdots x_n) &= Q(e, x_n) Q(x_n, x_{n-1}) \cdots Q(x_1, e), \quad  e \, x_1 \cdots x_n \in \hat S \\
	\hat Q(a, a^-) &= 1, \quad a \in \hat S - \{e\}
\end{align*}
It's easy to see $\hat Q$ is a valid transition probability function and corresponds to a downward run chain on $\hat S$.  Moreover, starting at $e$, we can define the two chains on a common probability space.  As before, if $\bs{\hat Y}$ is a process on $\hat S$, define the process $\bs{Y}$ on $S$ by $Y_n = m(\hat Y_n)$. The proof of the following proposition is straightforward.

\begin{proposition}
If $\bs{\hat Y}$ is a downward run chain on $\hat S$ with transition probability function $\hat Q$ and starting at $e$, then $\bs{Y}$ is a downward run chain on $S$ with transition probability function $Q$ and starting at $e$.
\end{proposition}

Of course the results of Section \ref{ss.downTrees} apply to $\bs{\hat Y}$. Also, with the chains coupled as in Proposition \ref{p.upCouple}, note that $T_e = \hat T_e$. The chains $\bs{Y}$ and $\bs{\hat Y}$ have the same classification.

Since $(\hat S, \preceq)$ is a rooted tree, $\hat G$ is an UPF.  On the other hand $G$ may not be an UPF on $S$, and conversely, there may exist UPFs $G$ on $S$ that cannot be realized as the standard invariant function for a recurrent downward run chain on $S$.

\section{Time Reversal}

The class of recurrent upward run chains and the class of downward run chains are time reversals of each other.  

\subsection{Reversing an upward run chain}

Suppose that $\bs{X}$ is a recurrent upward run chain with the transition probability function $P$ and standard invariant function $F$ (and other notation) as in Section \ref{s.upward}. Then the transition probability function $Q$ associated with the time reversed chain satisfies
\[ F(y)Q(y, x) = F(x) P(x, y), \quad x, \, y \in S \]
Thus, $Q(y, x) > 0$ if and only if $y = e$, or $y \ne e$ and $x \in B_y$. Hence the time reversed chain $\bs{Y}$ is a downward run chain and
\begin{equation}\label{eq.tr}
Q(y, x) = \frac{F(x)}{F(y)} P(x, y)
\end{equation}
Of course $F$ is also invariant for $Q$ and satisfies $F(e) = 1$ so $F = G$. In the positive recurrent case, $f = g$ and $\mu = \nu$.

We can interpret (\ref{eq.tr}) as
\[ Q(y, x) = \Pr_e[X_{T(y) - 1} = x | T_y \le T_e] \]
In particular, when $y = e$, we have
\[ Q(e, x) = \Pr_e[X_{T(e) - 1} = x] \]
Note that if $P(x, e)$ is constant in $x$, then the chains are positive recurrent, $f(x) = Q(e, x)$ is the invariant PDF.  That is,
\[ Q(e, x) = \frac{F(x)}{\mu(e)}, \quad x \in S \]

\subsection{Reversing a downward run chain}

Suppose that $\bs{Y}$ is a downward run chain with transition probability function $Q$ and standard invariant function $G$ (and other notation) as in Section \ref{s.downward}. The transition probability function $P$ associated with the time reversal satisfies
\[ G(x) P(x, y) = G(y) Q(y, x), \quad x, \, y \in S \]
Thus $P(x, y) > 0$ if and only if $y \in A_x$ or $y = e$. Hence, the time reversed chain $\bs{X}$ is an upward run chain and
\begin{equation} \label{eq.tr2}
P(x, y) = \frac{G(y)}{G(x)} Q(y, x)
\end{equation}
As before, $F = G$, and in the positive recurrent case, $f = g$, $\mu = \nu$.

We can interpret (\ref{eq.tr2}) as
\[ P(x, y) = \Pr_e[Y_{U(x) - 1} = y | U_x \le U_e] \]
When $y = e$, (\ref{eq.tr2}) gives
\begin{equation} \label{eq.tr3}
P(x, e) = \frac{Q(e, x)}{G(x)}
\end{equation}
which we can interpret as 
\[ P(x, e) = \Pr_e[U_x = 1 | U_x \le U_e] = \Pr_e(Y_1 = x | U_x \le U_e) \]
If $P(x, e)$ is constant in $x \in S$ then from (\ref{eq.tr3}), the chains are positive recurrent and $g(x) = Q(e, x)$ is the invariant PDF.  

\section{Examples and Special Cases} 

\subsection{Uniform Posets} \label{ss.uniform}

Suppose that the poset $(S, \preceq)$ is uniform, as defined in Section \ref{ss.uniform0}.  The general results in Sections \ref{s.upward} and \ref{s.downward} simplify significantly.

For the upward run chain $\bs{X} = (X_0, X_1, \ldots)$, note that
\[ F(x) = P^{d(e, x)}(e, x), \quad x \in S \]
That is, $F(x)$ is the probability, starting at $e$, that the chain moves strictly upward in $S$, reaching state $x$ in the minimum time $d(e, x)$. Thus,
\begin{align}
&F(x) = P^n(e, x), \quad x \in S_n, \, n \in \N \notag \\
&\Pr_e(T_e > n) = P^n(e, S_n), \quad n \in \N \label{eq.uni2} \\
&\mu(e) = \E_e(T_e) = \sum_{n=0}^\infty P^n(e, S_n) \label{eq.uni3}
\end{align}
The chain is recurrent if and only if the expression in (\ref{eq.uni2}) has limit 0 as $n \to \infty$ and positive recurrent if and only if the sum in (\ref{eq.uni3}) is finite. In the positive recurrent case, the invariant PDF $f$ is given by
\[ f(x) = \frac{P^n(e, x)}{\mu(e)}, \quad x \in S_n, \, n \in \N \]

Consider the downward run chain $\bs{Y} = (Y_0, Y_1, \cdots)$. If $Y_0 = e$ then $U_e = n + 1$ if and only if $Y_1 \in S_n$. Hence
\[ \Pr_e (U_e = n + 1) = \Pr_e(Y_1 \in S_n) = Q(e, S_n) \]
and so 
\begin{equation} \label{eq.nu}
\nu(e) := \E_e(U_e) = \sum_{n = 0}^\infty (n + 1) Q(x, S_n)
\end{equation}
Thus the chain is positive recurrent if and only if the sum in (\ref{eq.nu}) is finite. The standard invariant function $G$ also simplifies
\[ G(x) = \sum_{y \succeq x} Q(e, y) Q^{d(x, y)}(y, x), \quad x \in S \]

Now, for the upward run chain $\bs{X}$, consider the special case where 
\[ P(x, A_x) = \alpha_n, \; P(x, e) = 1 - \alpha_n; \quad x \in S_n, \, n \in \N \]
where $\alpha_n \in (0, 1)$ for $n \in \N$. Thus, the chain moves from level $n$ to level $n + 1$ with probability $\alpha_n$, and resets to $e$ with probability $1 - \alpha_n$. Let $P_+(x, y) = P(x, y) / \alpha_n$ for $x \in S_n$ and let $N_n = d(e, X_n)$ for $n \in \N$.  Then clearly $\bs{N} = (N_0, N_2, \ldots)$ is an ordinary success-runs chain on $(\N, +)$ while $P_+$ is the transition probability for a Markov chain that moves strictly upward in the poset. The behavior of the upward run chain $\bs{X}$ and the corresponding reversed downward run chain $\bs{Y}$ can be explained simply in terms of $\bs{N}$ and $P_+$.

Of course, the transition matrix $\hat P$ for $\bs{N}$ is given by 
\begin{equation*}
\hat P(n, n+1) = \alpha_n, \; \hat P(n, 0) = 1 - \alpha_n, \quad n \in \N
\end{equation*}
Clearly $\hat T_0 = T_e$ where $\hat T_0$ is the hitting time to 0 for the chain $\bs{N}$. The standard invariant function $\hat F$ for $\bs{N}$ is
\[ \hat F(n) = \alpha_0 \cdots \alpha_{n-1}, \quad n \in \N \]
Both chains are recurrent if and only if $\prod_{k=0}^\infty \alpha_k = 0$ and both chains are positive recurrent if and only if $\mu(0) = \sum_{n=0}^\infty \alpha_0 \cdots \alpha_{n-1} < \infty$. In the positive recurrent case, the invariant PDF for $\bs{N}$ is 
\[ \hat f(n) = \frac{\hat F (n)}{\mu(0)}, \quad n \in \N \]

Returning to the upward run chain $\bs{X}$, note that the standard invariant function $F$ satisfies
\[ F(x) = \hat F(n) P_+^n(e, x), \quad x \in S_n, \, n \in \N \]
In the positive recurrent case, the invariant PDF similarly satisfies
\[ f(x) = \hat f(n) P_+^n(e, x), \quad x \in S_n, \, n \in \N \]

The downward run chain obtained by reversing $\bs{N}$ has transition probabilities
\[ \hat Q(n+1, n) = 1, \; \hat Q(0, n) = \alpha_0 \cdots \alpha_{n-1} (1 - \alpha_n), \quad n \in \N \]
The downward run chain obtained by reversing $\bs{X}$ has transition probabilities
\begin{align*}
Q(y, x) &= \frac{P_+^n(e, x) P_+(x, y)}{P_+^{n+1}(e, y)}, \quad x \in S_n, y \in A_x, n \in \N \\
Q(e, x) &= \hat Q(0, n) P_+^n(e, x), \quad x \in S_n, n \in \N
\end{align*}
Note in particular that for $y \succ e$ and $x \in B_y$ the downward probability $Q(y, x)$ is independent of the parameters $(\alpha_0, \alpha_1, \ldots)$.

\subsection{Positive semigroups} \label{ss.ps}

Suppose now that the poset $(S, \preceq)$ is associated with a positive semigroup $(S, \cdot)$ with $I$ as the set of irreducible elements. It's natural to consider upward and downward runs that take advantage of the self-similarity noted in Section \ref{ss.semigroups0}.

An upward run chain $\bs{X}$ on $(S, \cdot)$ is {\em spatially homogeneous} if $P(x, xi) = r_i$, for all $x \in S$ and $i \in I$ where $r_i > 0$ and $r := \sum_{i \in I} r_i < 1$. It follows of course that $P(x, e) = 1 - r$ for $x \in S$. The standard invariant function for the corresponding chain on the free semigroup $(I^*, \cdot)$ is 
\[ F^*(i_1 \, i_2 \cdots i_n) = \prod_{k=1}^n r_{i_k}, \quad (i_1, i_2, \ldots, i_n) \in I^n \]
or equivalently $\ F^*(a) = \prod_{i \in I} r_i^{n_i(a)}$ where $n_i(a)$ is the number of times that letter $i$ appears in word $a$. The chain is positive recurrent and moreover, the invariant distribution is ``exponential'' (see \cite{Siegrist}). The corresponding standard invariant function $F$ for $\bs{X}$ is given by 
\[ F(x) = \sum_{i_1 i_2 \cdots i_n = x} F^*(i_1 i_2 \cdots i_n) \]
Of course, these results are a special case of the construction in Section \ref{s.upward}, since $(I^*, \cdot)$ is isomorphic to the path space $(\hat S, \preceq)$. 

For the downward run chain $\bs{Y}$ obtained by reversing $\bs{X}$, the function $x \mapsto Q(e, x)$ is the invariant PDF $f$, since $P(x, e)$ is constant in $x$.

\subsection{Upward and downward runs on $(\N^k, +)$}
Consider the case of the uniform positive semigroup $(\N^k, +)$. For $i = 1, \ldots, k$, let $\bs{u}_i \in \N^k$ be the element with 1 in position $i$ and 0 in all other positions; these are the irreducible elements.  For $\bs{x} \in \N^k$, let 
\[ C(\bs{x}) = \frac{\left(\sum_{i=1}^k x_i\right)!}{\prod_{i=1}^k x_i !} \]
This is a multinomial coefficient and gives the number of factorings of $\bs{x}$ over $I$; in each factoring, $\bs{u}_i$ must occur $x_i$ times.

Consider the upward run chain with uniform probabilities, as in Sections \ref{ss.uniform} and \ref{ss.ps}.  Specifically, let $r_i = P(\bs{x}, \bs{x} + \bs{u}_i)$, independent of $\bs{x} \in \N^k$, where $r := \sum_{i=1}^k p_i < 1$. For the upward run chain, the standard invariant function $F$ is given by
\[ F(\bs{x}) = C(\bs{x}) \prod_{i=1}^k r_i^{x_i}, \quad \bs{x} \in \N^k \]
The chain is positive recurrent and the invariant PDF is given by
\[ f(x) = (1 - r) F(\bs{x}) = (1 - r) C(\bs{x}) \prod_{i=1}^k r_i^{x_i}, \quad \bs{x} \in \N^k \]
If $\bs{Z} = (Z_1, \ldots, Z_k$) is a random vector with the invariant distribution then the (marginal) distribution of $Z_i$ is geometric on $\N$ with rate parameter 
\[ \alpha_i := 1 - \frac{r_i}{\sum\{r_j : j \ne i\}} \]

For the downward run chain $\bs{Y}$ obtained by reversing $\bs{X}$, the downward probabilities are given by 
\[ Q(\bs{x}, \bs{x} - \bs{u}_i) = \frac{x_i}{\sum_{j=1}^k x_j}, \quad \bs{x} \in \N - \{\bs{0}\} \]
Note that these probabilities are independent of $(r_1, r_2, \ldots, r_k)$.

\end{document}